\documentclass[10pt,onecolumn]{IEEEtran}
\usepackage[LGR,T1]{fontenc}
\usepackage[latin9]{inputenc}
\usepackage{color}
\usepackage{url}
\usepackage{amsmath}
\usepackage{amsthm}
\usepackage{amssymb}
\usepackage{graphicx}
\usepackage[unicode=true,
 bookmarks=true,bookmarksnumbered=true,bookmarksopen=true,bookmarksopenlevel=3,
 breaklinks=false,pdfborder={0 0 0},pdfborderstyle={},backref=false,colorlinks=true]
 {hyperref}
\hypersetup{
 pdfauthor={Lei Yu},
 pdfborderstyle={},pdfborderstyle={},pdfborderstyle={},pdfborderstyle={},pdfborderstyle={},pdfborderstyle={},pdfborderstyle={},pdfborderstyle={},pdfborderstyle={},pdfborderstyle={},pdfborderstyle={},pdfborderstyle={},pdfborderstyle={},pdfborderstyle={},pdfborderstyle={},pdfborderstyle={},pdfborderstyle={},pdfborderstyle={},pdfborderstyle={},pdfpagelayout=OneColumn,pdfnewwindow=true,pdfstartview=XYZ,plainpages=false,linkcolor=blue,urlcolor=blue,citecolor=red,anchorcolor=blue,linkcolor=blue,urlcolor=blue,citecolor=red,anchorcolor=blue}

\makeatletter

\newcommand{\lyxmathsym}[1]{\ifmmode\begingroup\def\b@ld{bold}
  \text{\ifx\math@version\b@ld\bfseries\fi#1}\endgroup\else#1\fi}

\theoremstyle{plain}
\newtheorem{thm}{\protect\theoremname}
\theoremstyle{definition}
\newtheorem{problem}[]{\protect\problemname}
\theoremstyle{plain}
\newtheorem{prop}[]{\protect\propositionname}
\theoremstyle{remark}
\newtheorem{rem}[]{\protect\remarkname}
\theoremstyle{plain}
\newtheorem{cor}[]{\protect\corollaryname}
\theoremstyle{definition}
\newtheorem{defn}[]{\protect\definitionname}
\theoremstyle{plain}
\newtheorem{lem}[]{\protect\lemmaname}

\usepackage{color}
\usepackage{cite}

  \providecommand{\definitionname}{Definition}
  \providecommand{\lemmaname}{Lemma}
  \providecommand{\remarkname}{Remark}
\providecommand{\theoremname}{Theorem}

\allowdisplaybreaks[1]
\providecommand{\corollaryname}{Corollary}
\providecommand{\lemmaname}{Lemma}
\providecommand{\propositionname}{Proposition}
\providecommand{\remarkname}{Remark}
\providecommand{\theoremname}{Theorem}

\providecommand{\corollaryname}{Corollary}
\providecommand{\lemmaname}{Lemma}
\providecommand{\propositionname}{Proposition}
\providecommand{\remarkname}{Remark}
\providecommand{\theoremname}{Theorem}

\makeatother

\providecommand{\corollaryname}{Corollary}
\providecommand{\definitionname}{Definition}
\providecommand{\lemmaname}{Lemma}
\providecommand{\problemname}{Problem}
\providecommand{\propositionname}{Proposition}
\providecommand{\remarkname}{Remark}
\providecommand{\theoremname}{Theorem}

\begin{document}
\title{An Improved Linear Programming Bound on the Average Distance of a
Binary Code}
\author{Lei Yu and Vincent Y. F. Tan\thanks{The authors are with the Department of Electrical and Computer Engineering,
National University of Singapore (Emails: \protect\url{leiyu@nus.edu.sg},
\protect\url{vtan@nus.edu.sg}).
V.~Y.~F. Tan is also with the Department of Mathematics, National
University of Singapore.} }
\maketitle
\begin{abstract}
Ahlswede and Katona (1977) posed the following  isodiametric problem
in Hamming spaces: For every $n$ and $1\le M\le2^{n}$, determine
the minimum average Hamming distance of binary codes with length $n$
and size $M$.  Fu, Wei, and Yeung (2001) used linear programming
duality to derive a lower bound on the minimum average distance. However,
their linear programming approach was not completely exploited. In
this paper, we improve Fu-Wei-Yeung's bound by finding a better feasible
solution to their dual program. For fixed $0<a\le1/2$ and for $M=\left\lceil a2^{n}\right\rceil $,
our feasible solution attains the asymptotically optimal value of
Fu-Wei-Yeung's dual program as $n\to\infty$. Hence for $0<a\le1/2$,
all possible asymptotic bounds that can be derived by Fu-Wei-Yeung's
linear program have been characterized. Furthermore, noting that the
average distance of a code is closely related to weights of Fourier
coefficients of a Boolean function, we also apply the linear programming
technique to prove bounds on Fourier weights of a Boolean function
of various degrees.
\end{abstract}


\section{\label{sec:Introduction}Introduction}

A binary $(n,M)$-code is a subset $A$ of $\{-1,1\}^{n}$ with size
$M$. The average distance of $A$ is defined to be the average Hamming
distance of every pair of codewords in $A$. Ahlswede and Katona \cite{ahlswede1977contributions}
posed the following problem concerning the extremal combinatorics
in Hamming space: For every $1\le M\le2^{n}$, determine the minimum
of the average distance $D\left(A\right)$ over all sets $A\subseteq\{-1,1\}^{n}$
of a given cardinality $M$. K\"undgen \cite{kundgen1998covering} observed
that this problem is equivalent to a covering problem in graph theory.
Ahlswede and Alth\"ofer \cite{ahlswede1994asymptotic} considered the
case in which the size of code increases exponentially in $n$ and
the exponent is strictly between $0$ and $1$. They provided nearly
optimal solutions (which are attained by Hamming spheres) to Ahlswede-Katona's
problem for the asymptotic case in which $n\to\infty$.  Using a
linear programming approach, Mounits \cite{mounits2008lower} studied
codes whose sizes are linear in $n$ (i.e., codes with ``small''
sizes). He showed that when the size of code is $2n$, the asymptotic
value of the minimum average distance is $\frac{5}{2}$ as $n\to\infty$.
Alth\"ofer and Sillke \cite{althofer1992average}, Fu, Xia, together
with other authors \cite{shutao1998average,fu1997expectation,fu1999hamming,fu2001minimum},
as well as Mounits \cite{mounits2008lower}, proved various bounds
on the minimum average distance, which are sharp in certain regimes
when the code size is ``large'' (e.g., $M=2^{n-1}$ or $2^{n-2}$).
In particular, Fu, Wei, and Yeung \cite{fu2001minimum} used linear
programming duality to show that for any $(n,M)$-code $A$ such that
$a:=\frac{M}{2^{n}}\le\frac{1}{2}$, 
\begin{equation}
D\left(A\right)\geq\frac{n}{2}-\frac{1}{4a},\label{eq:-Fu}
\end{equation}
and equality in \eqref{eq:-Fu} holds for $M=2^{n-1}$ or $2^{n-2}$
by setting $A$ to be a subcube (e.g., $A=\{1\}\times\left\{ -1,1\right\} ^{n-1}$
for $M=2^{n-1}$ and $A=\{1\}^{2}\times\left\{ -1,1\right\} ^{n-2}$
for $M=2^{n-2}$). In Fu-Wei-Yeung's linear programming approach,
it was observed that minimizing the average distance over all $(n,M)$-codes
is equivalent to minimizing the average distance over all possible
dual distance distributions of $(n,M)$-codes. By relaxing the condition
that the dual distance distribution lies in a certain finite subset
of the nonnegative orthant $\mathbb{R}_{\ge0}^{n+1}$ to the condition
that it can be any vector in $\mathbb{R}_{\ge0}^{n+1}$, the latter
minimization problem is shown to be equivalent to a linear program.
By strong duality of linear programming, the optimal value of this
linear program is equal to that of its dual problem. On the other
hand, the optimal value of the dual (maximization) program can be
lower bounded by evaluating the dual objective at a feasible solution.
This results in a lower bound for the original problem (i.e., the
minimum average distance problem). Moreover, a better feasible solution
will result in a tighter bound for the original problem. Hence finding
a good solution to the dual program is particularly important in this
approach. In \cite{fu2001minimum}, Fu, Wei, and Yeung derived the
bound \eqref{eq:-Fu} by finding a simple feasible solution $\left(0,0,...,0,\frac{1}{2}\right)$.
(Note that this feasible solution is independent of the parameter
$a$.) In this paper, we improve Fu-Wei-Yeung's bound. We first find
a better feasible solution to the dual program, and then prove that
our feasible solution is asymptotically optimal as $n\to\infty$.
Hence all possible bounds that can be derived by using Fu-Wei-Yeung's
linear programming approach are characterized asymptotically.\footnote{Note that we are \emph{not} the first to study the asymptotic optimality
of a specific linear programming approach. In coding theory, McEliece,
Rodemich, Rumsey, and Welch \cite{mceliece1977new} provided the best
known upper bound for the sphere packing problem in Hamming spaces.
This bound was obtained by finding a feasible solution to the dual
program in Delsarte's linear programming approach. In \cite{samorodnitsky2001optimum},
Samorodnitsky studied the optimality of Delsarte's linear programming
approach, and conjectured that McEliece-Rodemich-Rumsey-Welch's feasible
solution is an asymptotically optimal solution to the dual program
in Delsarte's linear programming approach as the blocklength $n\to\infty$.}

The average distance of a code is closely related to the (Fourier)
weight of a Boolean function at degree $1$. For a Boolean function
$f:\{-1,1\}^{n}\to\{-1,1\}$, we use $\hat{f}_{S}$ for\footnote{Throughout this paper, we denote $[m:n]:=\left\{ m,m+1,...,n\right\} $.}
$S\subseteq[1:n]$ to denote Fourier coefficients of $f$. Then the
degree-$1$ Fourier weight of $f$ is defined as 
\begin{align*}
\mathbf{W}_{1} & :=\sum_{S:|S|=1}\hat{f}_{S}^{2}.
\end{align*}
The degree-$1$ Fourier weight $\mathbf{W}_{1}$ and the average distance
of $A=f^{-1}(1):=\left\{ \mathbf{x}\in\{-1,1\}^{n}:f\left(\mathbf{x}\right)=1\right\} $
admit the following intimate relationship \cite{yu2019bounds}:
\begin{align}
\mathbf{W}_{1} & =4a^{2}\left(n-2D\left(A\right)\right).\label{eq:-8}
\end{align}
Hence the estimation of the degree-$1$ Fourier weight of $f$ is
equivalent to the estimation of the average distance of $A$. It is
worth noting that the estimation of Fourier coefficients of a Boolean
function is an important topic in theoretical computer science and
Fourier analysis, which has found many applications in coding theory,
noise-sensitivity theory, and combinatorics \cite{O'Donnell14analysisof,green2008boolean,chang2002polynomial,friedgut2002boolean,defant2019fourier}.
 In this paper, we also apply the linear programming technique to
prove upper bounds on the degree-$m$ Fourier weight of a Boolean
function for different $m$'s.

This paper is organized as follows. In Section II, we introduce some
background concerning the minimum average distance problem. Specifically,
we provide the definitions of several quantities (including the distance
distribution, the average distance, and the distance enumerator) and
briefly describe Fu-Wei-Yeung's linear programming approach in \cite{fu2001minimum}.
In Section III, we improve Fu-Wei-Yeung's bound by finding a new feasible
solution to their dual program. The asymptotic optimality of our feasible
solution is also studied. Furthermore, we also compare our improved
linear programming bound with existing bounds, including Chang's bound
\cite[Lemma 3.1]{chang2002polynomial} and the hypercontractivity
bound \cite{yu2019bounds}. In Section IV, we apply linear programming
approach to obtain upper bounds on the degree-$m$ Fourier weight
of a Boolean function. Finally, in Section V, we apply  our results
to  estimate the noise stability of Boolean functions.

\section{Background}

\subsection{Definitions}

For a subset of the Boolean hypercube (termed a \emph{code}) $A\subseteq\{-1,1\}^{n}$,
the \emph{distance distribution} of $A$ is the following probability
mass function: 
\[
P^{\left(A\right)}(i):=\frac{1}{|A|^{2}}\left|\left\{ \left(\mathbf{x},\mathbf{x}'\right)\in A^{2}:d_{\mathrm{H}}\left(\mathbf{x},\mathbf{x}'\right)=i\right\} \right|,\;i\in[0:n],
\]
where $d_{\mathrm{H}}\left(\mathbf{x},\mathbf{x}'\right):=\left|\left\{ i:\:x_{i}\neq x'_{i}\right\} \right|$
denotes the Hamming distance between vectors $\mathbf{x}$ and $\mathbf{x}'$
(i.e., the number of components of $\mathbf{x}$ and $\mathbf{x}'$
that differ). It is clear that $P^{\left(A\right)}(0)=\frac{1}{|A|}$,
$\sum_{i=0}^{n}P^{\left(A\right)}(i)=1$, and $P^{\left(A\right)}(i)\ge0$
for $i\in[0:n]$.

Define the \emph{average distance }of the code $A\subseteq\{-1,1\}^{n}$
as 
\[
D\left(A\right):=\frac{1}{|A|^{2}}\sum_{\left(\mathbf{x},\mathbf{x}'\right)\in A^{2}}d_{\mathrm{H}}\left(\mathbf{x},\mathbf{x}'\right)=\sum_{i=0}^{n}P^{\left(A\right)}(i)\cdot i.
\]

Define the \emph{distance enumerator} of $A\subseteq\{-1,1\}^{n}$,
with $z$ as the indeterminate, as 
\[
\Gamma_{z}\left(A\right):=\frac{1}{|A|^{2}}\sum_{\left(\mathbf{x},\mathbf{x}'\right)\in A^{2}}z^{d_{\mathrm{H}}\left(\mathbf{x},\mathbf{x}'\right)}=\sum_{i=0}^{n}P^{\left(A\right)}(i)\cdot z^{i}.
\]
Clearly, $\Gamma_{z}\left(A\right)$ is the the generating function
of $P^{\left(A\right)}$. For $z=1$, $\Gamma_{1}\left(A\right)=1$.

The \emph{dual distance distribution} of $A$ is defined by 
\begin{equation}
Q^{\left(A\right)}(i):=\frac{1}{|A|^{2}}\sum_{\mathbf{u}\in\{0,1\}^{n}:w_{\mathrm{H}}(\mathbf{u})=i}\left(\sum_{\mathbf{x}\in\frac{A+1}{2}}\left(-1\right)^{\left\langle \mathbf{u},\mathbf{x}\right\rangle }\right)^{2},\;i\in[0:n],\label{eq:-30}
\end{equation}
where $w_{\mathrm{H}}(\mathbf{u}):=d_{\mathrm{H}}\left(\mathbf{u},\mathbf{0}\right)$
denotes the Hamming weight (i.e., the number of nonzero components)
of a vector $\mathbf{u}$, and $\left\langle \mathbf{u},\mathbf{x}\right\rangle :=\left(\sum_{i=1}^{n}u_{i}x_{i}\right)\,\mathrm{mod}\,2$
denotes the inner product of vectors $\mathbf{u}$ and $\mathbf{x}$
in $\mathbb{F}_{2}^{n}$. Clearly, 
\begin{align}
Q^{\left(A\right)}(0) & =1,\label{eq:-58}\\
Q^{\left(A\right)}(i) & \ge0\textrm{ for }i\in[0:n].\label{eq:-59}
\end{align}

The \emph{dual distance enumerator} of $A$ is defined as 
\begin{equation}
\Pi_{z}\left(A\right):=\sum_{i=0}^{n}Q^{\left(A\right)}(i)\cdot z^{i},\quad z\ge0.\label{eq:-51}
\end{equation}
The following MacWilliams\textendash Delsarte identities hold \cite{macwilliams1977theory}.
\begin{align}
\Pi_{z}\left(A\right) & =\left(1+z\right)^{n}\Gamma_{\frac{1-z}{1+z}}\left(A\right)\label{eq:-22}\\
\Gamma_{z}\left(A\right) & =\left(\frac{1+z}{2}\right)^{n}\Pi_{\frac{1-z}{1+z}}\left(A\right).\label{eq:-23}
\end{align}
By \eqref{eq:-22}, 
\begin{equation}
\sum_{i=0}^{n}Q^{\left(A\right)}(i)=\frac{2^{n}}{|A|}.\label{eq:-26}
\end{equation}
 Hence for this case, $\frac{|A|}{2^{n}}Q^{\left(A\right)}(\cdot)$
is a probability mass function.

Consider the Fourier basis $\left\{ \chi_{S}\right\} _{S\subseteq[1:n]}$
with $\chi_{S}(\mathbf{x}):=\prod_{i\in S}x_{i}$ for $S\subseteq[1:n]$.
Then for a Boolean function $f:\{-1,1\}^{n}\to\{-1,1\}$, define its
Fourier coefficients as 
\begin{align}
 & \hat{f}_{S}:=\mathbb{E}_{\mathbf{x}\sim\mathrm{Unif}\left\{ -1,1\right\} ^{n}}[f(\mathbf{x})\chi_{S}(\mathbf{x})],\;S\subseteq[1:n].\label{eq:-33}
\end{align}
Then the Fourier expansion of a Boolean function $f$ (cf. \cite[Equation (1.6)]{O'Donnell14analysisof})
is
\begin{align*}
 & f(\mathbf{x})=\sum_{S\subseteq[1:n]}\hat{f}_{S}\chi_{S}(\mathbf{x}).
\end{align*}
The \emph{degree-$m$ Fourier weight }of $f$ is defined as 
\begin{align*}
\mathbf{W}_{m} & :=\sum_{S:|S|=m}\hat{f}_{S}^{2},\quad m\in[0:n].
\end{align*}
By definition, it is easily seen that
\begin{align*}
\sum_{m=0}^{n}\mathbf{W}_{m} & =1,\\
\mathbf{W}_{0} & =\left(2a-1\right)^{2}
\end{align*}
where $a=\left|f^{-1}(1)\right|/2^{n}$.

For a code $A\subseteq\{-1,1\}^{n}$,  the dual distribution of $A$
and the Fourier coefficients of $f=2\cdot1_{A}-1$ admit the following
relationship \cite{yu2019bounds}:
\begin{equation}
Q^{\left(A\right)}(k)=\begin{cases}
1, & k=0\\
\frac{1}{4a^{2}}\mathbf{W}_{k}, & 1\le k\le n
\end{cases},\label{eq:-52}
\end{equation}
where $a=\frac{|A|}{2^{n}}$. For $k=1$, 
\begin{align}
\mathbf{W}_{1} & =4a^{2}Q^{\left(A\right)}(1)\label{eq:-32}\\
 & =4a^{2}\left(n-2D\left(A\right)\right).\label{eq:-27}
\end{align}

\subsection{\label{subsec:Fu-Wei-Yeung's-Linear-Programmin}Fu-Wei-Yeung's Linear
Programming Approach}

For each $k\in[0:n]$ and indeterminate $x$, the Krawtchouk polynomials
\cite{macwilliams1977theory} are defined as\footnote{Here the (generalized) binomial coefficients ${x \choose j}:=\frac{x(x-1)\cdots(x-j+1)}{j!}$.}
\[
K_{k}^{(n)}\left(x\right):=\sum_{j=0}^{k}(-1)^{j}{x \choose j}{n-x \choose k-j},
\]
whose generating function satisfies 
\begin{equation}
\sum_{k=0}^{\infty}K_{k}^{(n)}\left(x\right)z^{k}=(1-z)^{x}(1+z)^{n-x}.\label{eq:-22-1}
\end{equation}
For brevity and if there is no ambiguity, we denote $K_{k}^{(n)}$
as $K_{k}$.

For $i=0,1$, we have 
\begin{align*}
K_{k}\left(0\right) & ={n \choose k}
\end{align*}
and 
\begin{align*}
K_{k}\left(1\right) & ={n \choose k}\left(1-\frac{2k}{n}\right).
\end{align*}
Combining \eqref{eq:-22}, \eqref{eq:-23}, and \eqref{eq:-22-1}
yields that the distance distribution and its dual are related via
the Krawtchouk transform as shown in the following:
\begin{align}
Q^{\left(A\right)}(k) & =\sum_{i=0}^{n}P^{\left(A\right)}(i)K_{k}\left(i\right)\label{eq:-22-2}\\
P^{\left(A\right)}(k) & =\frac{1}{2^{n}}\sum_{i=0}^{n}Q^{\left(A\right)}(i)K_{k}\left(i\right).\label{eq:-23-1}
\end{align}

Given a code $A$ of size $M$, by \eqref{eq:-26} and \eqref{eq:-27},
the average distance of $A$ satisfies \cite[Section 4.1]{fu2001minimum}
\begin{align}
D\left(A\right) & =\frac{n+1}{2}-\frac{1}{2a}+\frac{1}{2}\sum_{i=2}^{n}Q^{\left(A\right)}(i),\label{eq:-71}
\end{align}
where $a=\frac{|A|}{2^{n}}$. Hence minimizing $D\left(A\right)$
is equivalent to minimizing $\sum_{i=2}^{n}Q^{\left(A\right)}(i)$.
Recall that $Q^{\left(A\right)}(\cdot)$ denotes the dual distance
distribution of $A$, which satisfies \eqref{eq:-58} and \eqref{eq:-59}.
By \eqref{eq:-23-1},
\begin{equation}
\sum_{i=0}^{n}Q^{\left(A\right)}(i)K_{k}\left(i\right)\ge0.\label{eq:-77}
\end{equation}
In \cite{fu2001minimum}, Fu, Wei, and Yeung considered a relaxed
version of the minimization (integer program) of $\sum_{i=2}^{n}Q^{\left(A\right)}(i)$
over the dual distance distribution $Q^{\left(A\right)}$.  Instead
of the discrete optimization of $\sum_{i=2}^{n}Q^{\left(A\right)}(i)$
(since given $n$, there are only finitely many codes and the corresponding
dual distance distributions), they allowed $\left(Q^{\left(A\right)}(0),Q^{\left(A\right)}(1),...,Q^{\left(A\right)}(n)\right)$
to be any nonnegative vector $\left(u_{0},u_{1},...,u_{n}\right)$
such that
\begin{align*}
 & u_{0}=1,u_{i}\ge0,\;i\in\left[2:n\right];\\
 & \sum_{i=0}^{n}u_{i}=\frac{1}{a};\\
 & \sum_{i=0}^{n}u_{i}K_{k}\left(i\right)\ge0,\;k\in\left[0:n\right].
\end{align*}
Then in order to minimize $\sum_{i=2}^{n}Q^{\left(A\right)}(i)$,
they considered the following linear program.

\begin{problem}
\label{prob:Primal-Problem:}Primal Problem:
\[
\Lambda(n;a):=\min_{u_{2},u_{3},...,u_{n}}\sum_{i=2}^{n}u_{i}
\]
subject to the inequalities 
\begin{align*}
 & u_{i}\ge0,\;i\in\left[2:n\right];\\
 & \sum_{i=2}^{n}\left[K_{k}\left(1\right)-K_{k}\left(i\right)\right]u_{i}\le K_{k}\left(0\right)+K_{k}\left(1\right)\left(\frac{1}{a}-1\right),\;k\in\left[1:n\right].
\end{align*}
\end{problem}
The dual is the following optimization problem.
\begin{problem}
\label{prob:Dual-Problem:}Dual Problem:
\begin{equation}
\overline{\Lambda}(n;a):=\max_{x_{1},x_{2},...,x_{n}}-\sum_{k=1}^{n}\left[K_{k}\left(0\right)+K_{k}\left(1\right)\left(\frac{1}{a}-1\right)\right]x_{k}\label{eq:-55}
\end{equation}
subject to the inequalities 
\begin{align*}
 & x_{k}\ge0,\;k\in\left[1:n\right];\\
 & \sum_{k=1}^{n}\left[K_{k}\left(1\right)-K_{k}\left(i\right)\right]x_{k}\ge-1,\;i\in\left[2:n\right].
\end{align*}
\end{problem}
By strong duality in linear programming,\footnote{Obviously, in the primal problem, since $u_{i}\ge0$, the primal problem
is bounded. On the other hand, the existence of a code $A$ with size
$M:=a2^{n}$ ensures that $u_{i}=Q^{\left(A\right)}(i)$ is a feasible
solution. Hence the primal problem has an optimal solution.} $\Lambda(n;a)=\overline{\Lambda}(n;a)$. Using this linear programming
approach, Fu, Wei, and Yeung obtained the following important result.
\begin{thm}
\label{thm:-For-any}\cite{fu2001minimum} For any code $A$ of size
$M$, 
\[
\sum_{i=2}^{n}Q^{\left(A\right)}(i)\geq\overline{\Lambda}(n;a).
\]
\end{thm}
By \eqref{eq:-71} and Theorem \ref{thm:-For-any}, 
\begin{align}
D\left(A\right)-\frac{n}{2} & =\frac{1}{2}-\frac{1}{2a}+\frac{1}{2}\sum_{i=2}^{n}Q^{\left(A\right)}(i)\label{eq:-71-1}\\
 & \geq\frac{1}{2}-\frac{1}{2a}+\frac{1}{2}\overline{\Lambda}(n;a).\label{eq:-76}
\end{align}
In \cite{fu2001minimum}, Fu, Wei, and Yeung found a simple feasible
solution $\left(0,0,...,0,\frac{1}{2}\right)$ to Problem \ref{prob:Dual-Problem:}.
Substituting this feasible solution into the dual objective function
in \eqref{eq:-55}, they obtained the lower bound $\frac{1}{2a}-1$
on $\overline{\Lambda}(n;a)$. This solution leads to the lower bound
in \eqref{eq:-Fu} on the average distance.

\section{Improved Linear Programming Bounds}

In this section, we first improve Fu-Wei-Yeung's bound. We then compare
our new bound with several existing bounds.

\subsection{Improved Linear Programming Bounds}

It was shown in \cite{fu2001minimum} that 
\begin{equation}
|A|^{2}D\left(A\right)-|A^{c}|^{2}D\left(A^{c}\right)=\left(|A|-|A^{c}|\right)n2^{n-1}.\label{eq:-11-1-1}
\end{equation}
This implies that bounding $D\left(A\right)$ is equivalent to bounding
$D\left(A^{c}\right)$. Hence it suffices to consider code sizes $M\leq2^{n-1}$,
i.e., $a:=\frac{M}{2^{n}}\le\frac{1}{2}$. We next provide a simple
observation for the average distance. The proof of Proposition \ref{prop:monotonicity}
is provided in Appendix \ref{sec:Proof-of-Theorem-ImprovedLPB}.
\begin{prop}[Monotonicity]
\label{prop:monotonicity} For $a:=\frac{M}{2^{n}}\le\frac{1}{2}$
and $k\in\mathbb{Z}_{\ge0}$, 
\begin{equation}
\min_{A\subseteq\{-1,1\}^{n+k}:|A|=2^{k}M}D\left(A\right)-\frac{k}{2}\leq\min_{A\subseteq\{-1,1\}^{n}:|A|=M}D\left(A\right).\label{eq:-56-3}
\end{equation}
\end{prop}
By induction, this proposition implies that for a dyadic rational
$a=\frac{M}{2^{n}}$, the sequence $\min_{A\subseteq\{-1,1\}^{n+k}:|A|=2^{k}M}D\left(A\right)-\frac{k}{2}$
is non-increasing in $k$.

Now we turn to provide the promised improvement of Fu-Wei-Yeung's
bound. 
\begin{thm}[Improved LP Bound]
\label{thm:ImprovedLPB} For $a:=\frac{M}{2^{n}}\le\frac{1}{2}$,
\begin{equation}
\min_{A:|A|=M}D\left(A\right)\geq\frac{n}{2}-\varphi(a),\label{eq:-56}
\end{equation}
where 
\[
\varphi(a):=\begin{cases}
\frac{1}{\sqrt{a}}-1 & 0\le a\le\frac{1}{4}\\
\frac{1}{4a} & \frac{1}{4}<a\le\frac{1}{2}
\end{cases}.
\]
\end{thm}
The proof of Theorem \ref{thm:ImprovedLPB} is provided in Appendix
\ref{sec:Proof-of-Theorem-ImprovedLPB}. In this proof, we in fact
show that 
\begin{equation}
\liminf_{n\to\infty}\overline{\Lambda}(n;a)\geq\theta(a),\label{eq:-72}
\end{equation}
where 
\begin{equation}
\theta(a):=\begin{cases}
\frac{\left(1-\sqrt{a}\right)^{2}}{a} & 0\le a<1/4\\
\frac{1}{2a}-1 & 1/4\le a\le1/2
\end{cases}.\label{eq:theta}
\end{equation}
Combining \eqref{eq:-76} and \eqref{eq:-72} yields that 
\begin{align}
\liminf_{n\to\infty}\left\{ \min_{A:|A|=a2^{n}}D\left(A\right)-\frac{n}{2}\right\}  & \geq\frac{1}{2}-\frac{1}{2a}+\frac{1}{2}\liminf_{n\to\infty}\overline{\Lambda}(n;a)\label{eq:-76-1}\\
 & \geq\frac{1}{2}-\frac{1}{2a}+\frac{1}{2}\theta(a).
\end{align}
By Proposition \ref{prop:monotonicity}, \eqref{eq:-56} follows.

The bound for the case of $\frac{1}{4}<a\le\frac{1}{2}$ was proven
by Fu, Wei, and Yeung and stated in \eqref{eq:-Fu}. This bound was
proved by substituting the dual feasible solution $\left(0,0,...,0,\frac{1}{2}\right)$
into the dual objective function of Problem \ref{prob:Dual-Problem:}.
 In our proof, we constructed another feasible solution 
\begin{equation}
\mathbf{x}^{*}=\left(0,...,0,x_{k}^{*},x_{k+1}^{*},0,...,0\right)\label{eq:oursolution-2}
\end{equation}
 with
\begin{align}
x_{k}^{*} & =\frac{1+2\left(\frac{k}{n}+\frac{1}{n}\right)^{2}-2\left(\frac{k}{n}+\frac{1}{n}\right)-\frac{1}{n}}{{n \choose k}2\frac{k}{n}\left(2\frac{k}{n}-1+\frac{1}{n}\left(2\frac{k}{n}+2\frac{1}{n}-1\right)\right)}\label{eq:oursolution}\\
x_{k+1}^{*} & =\frac{1-\frac{k}{n}}{{n \choose k+1}\left(2\frac{k}{n}-1+\frac{1}{n}\left(2\frac{k}{n}+2\frac{1}{n}-1\right)\right)},\label{eq:oursolution-1}
\end{align}
where $k=2\left\lfloor \frac{\beta n}{2}\right\rfloor $ for some
$\beta\in\left(\frac{1}{2},1\right)$. Here the $\beta$ we chose
is 
\[
\beta=\begin{cases}
\frac{1}{2\left(1-\sqrt{a}\right)} & 0\le a<1/4\\
1 & 1/4\le a\le1/2
\end{cases}
\]
which depends on the value of $a$.  For fixed $\beta$, letting
$n\to\infty$, we have 
\begin{align*}
{n \choose k}x_{k}^{*} & \to\frac{1+2\beta^{2}-2\beta}{2\beta\left(2\beta-1\right)}\\
{n \choose k+1}x_{k+1}^{*} & \to\frac{1-\beta}{2\beta-1}.
\end{align*}
In our proof, we show that this sequence of feasible solutions $\mathbf{x}^{*}$
(indexed by $n$) leads to the bound \eqref{eq:-72}.

One may wonder whether it is possible to further improve the bound
in \eqref{eq:-56} by finding more complicated dual feasible solutions
(our solution is only $2$-sparse).  In the following, we show that
the answer is no. To show this, we first prove the following bounds
for Problem \ref{prob:Dual-Problem:}. The proof of Proposition \ref{thm:Optimality-2}
is provided in Appendix \ref{sec:Proof-of-Theorem-Optimality}.
\begin{prop}[Bounds on $\overline{\Lambda}(n;a)$]
\label{thm:Optimality-2} We have the following bounds on $\overline{\Lambda}(n;a)$,
defined in \eqref{eq:-55}.
\begin{enumerate}
\item For  $a\leq1/4$ and $n\ge1/a-1$, we have 
\begin{equation}
\overline{\Lambda}(n;a)\leq\frac{\left(1-\frac{1}{n}\right)\frac{s}{a}}{\frac{1+s}{1-a}\left(\frac{1+s}{1-a}-1-\frac{1}{n}\right)},\label{eq:-75}
\end{equation}
where $s:=\sqrt{a-\frac{1-a}{n}}$. 
\item For $a>1/4$ and $n\ge\frac{1-a}{a-\left(1-2a\right)^{2}}$, we have
\begin{align}
\overline{\Lambda}(n;a) & \leq\frac{1}{2a}-1.\label{eq:-74}
\end{align}
\end{enumerate}
\end{prop}
\begin{rem}
Fu, Wei, and Yeung showed that for any $a,n$, $\left(0,0,...,0,\frac{1}{2}\right)$
is a feasible solution to Problem \ref{prob:Dual-Problem:}. This
solution when substituted into the dual objective function yields
the value of $\frac{1}{2a}-1$. Hence combining this with \eqref{eq:-74},
we have that for $a>1/4$ and $n\ge\frac{1-a}{a-\left(1-2a\right)^{2}}$,
\begin{align}
\overline{\Lambda}(n;a) & =\frac{1}{2a}-1,\label{eq:-74-1}
\end{align}
and so for $a>1/4$, Fu-Wei-Yeung's solution is optimal for sufficiently
large $n$.
\end{rem}
Observe that the RHS of \eqref{eq:-75} satisfies that

\[
\frac{\left(1-\frac{1}{n}\right)\frac{s}{a}}{\frac{1+s}{1-a}\left(\frac{1+s}{1-a}-1-\frac{1}{n}\right)}\to\frac{\left(1-\sqrt{a}\right)^{2}}{a}
\]
as $n\to\infty$. Therefore,
\begin{equation}
\limsup_{n\to\infty}\overline{\Lambda}(n;a)\leq\theta(a),\label{eq:-72-2}
\end{equation}
where $\theta(a)$ is defined in \eqref{eq:theta}. Combining this
with \eqref{eq:-72}, we obtain the following theorem. 
\begin{thm}[Asymptotic Optimality of the Bound in \eqref{eq:-56}]
\label{thm:Optimality} For fixed $a\le\frac{1}{2}$ and for\footnote{Throughout this paper, we denote $\left\lceil x\right\rceil $ as
the least integer greater than or equal to $x$ and $\left\lfloor x\right\rfloor $
as the greatest integer smaller than or equal to $x$.} $M=\left\lceil a2^{n}\right\rceil $, 
\begin{equation}
\lim_{n\to\infty}\overline{\Lambda}(n;a)=\theta(a),\label{eq:-72-1}
\end{equation}
and the sequence of vectors $\left\{ \mathbf{x}^{*}\right\} $ defined
in \eqref{eq:oursolution-2}\textendash \eqref{eq:oursolution-1}
asymptotically attains $\theta(a)$ in \eqref{eq:-72-1}. 
\end{thm}
Recall the relationship between $D\left(A\right)$ and $\overline{\Lambda}(n;a)$
in \eqref{eq:-76}. By Theorem \ref{thm:Optimality}, the lower bound
in \eqref{eq:-76} satisfies that 
\begin{align*}
 & \lim_{n\to\infty}\left\{ \frac{1}{2}-\frac{1}{2a}+\frac{1}{2}\overline{\Lambda}(n;a)\right\} \\
 & =\frac{1}{2}-\frac{1}{2a}+\frac{1}{2}\theta(a)\\
 & =-\varphi(a).
\end{align*}
This means that for fixed $a\in(0,\frac{1}{2})$, Fu-Wei-Yeung's
linear programming approach cannot be used to obtain a bound that
is better than \eqref{eq:-56} asymptotically as $n\to\infty$. In
other words, our proposed sequence of $2$-sparse solutions $\left\{ \mathbf{x}^{*}\right\} $
is asymptotically optimal in terms of maximizing the dual objective
function in Problem \ref{prob:Dual-Problem:}.

\subsection{Comparisons to Other Bounds}

Chang proved the following bound by using results in additive combinatorics
\cite{chang2002polynomial,O'Donnell14analysisof}. Beautiful information-theoretic
proofs of the same bound were provided by Impagliazzo, Moore, and
Russell \cite{impagliazzo2014entropic,hambardzumyan2019chang} as
well as Hambardzumyan and Li \cite{hambardzumyan2019chang}.
\begin{prop}[{Chang\textquoteright s Bound \cite[Lemma 3.1]{chang2002polynomial}}]
\label{prop:Chang} For $1\le M\le2^{n}$ and $a=\frac{M}{2^{n}}$,
we have 
\begin{equation}
\min_{A:|A|=M}D\left(A\right)\geq\frac{n}{2}-\ln\frac{1}{a}.\label{eq:-31}
\end{equation}
\end{prop}
By using hypercontractivity inequalities, in a recent paper \cite{yu2019bounds}
the present authors showed the following bound on the average distance.

\begin{prop}[Hypercontractivity Bound \cite{yu2019bounds}]
\label{prop:RHB-1} For $1\le M\le2^{n}$, we have 
\begin{equation}
\min_{A:|A|=M}D\left(A\right)\geq\frac{n}{2}-\psi\left(a\right),\label{eq:-39}
\end{equation}
where 
\begin{align}
\psi\left(a\right) & :=\inf_{t>0,t\neq1}\frac{\left(ta+\overline{a}\right)\left[at\ln t-\left(ta+\overline{a}\right)\ln\left(ta+\overline{a}\right)\right]}{a^{2}\left(t-1\right)^{2}}.\label{eq:-45-1}
\end{align}
\end{prop}
As shown in \cite{yu2019bounds}, the hypercontractivity bound is
tighter than Chang's bound for all $a\in(0,1]$.

\begin{figure}
\centering \includegraphics[scale=0.7]{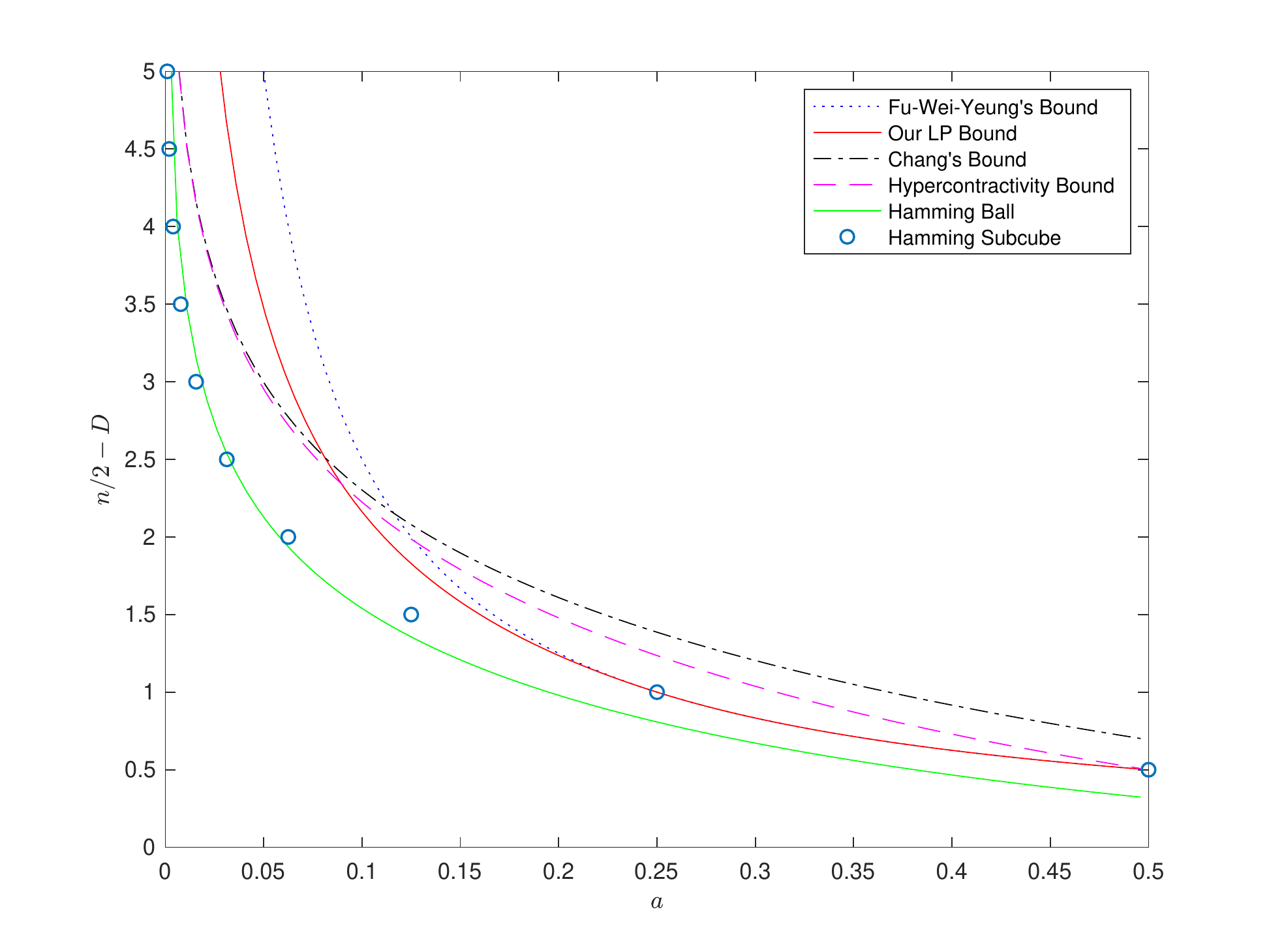}

\includegraphics[scale=0.7]{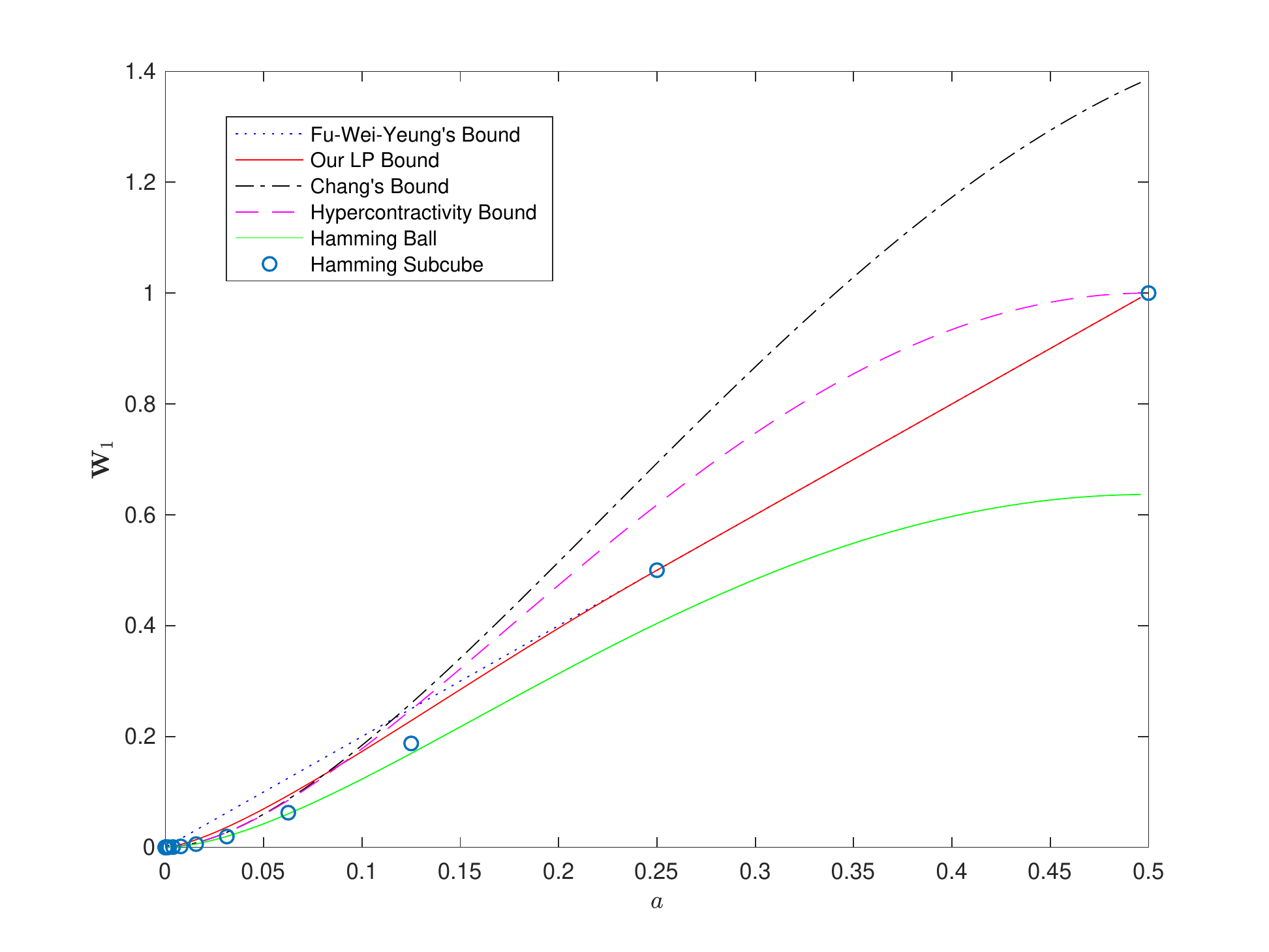}

\caption{\label{fig:simulation-1}Illustration of Fu-Wei-Yeung's bound in \eqref{eq:-Fu},
the improved linear programming bound in \eqref{eq:-56}, Chang's
bound in \eqref{eq:-31}, and the hypercontractivity bound in \eqref{eq:-39}.
In the top figure, the vertical axis corresponds to the gap $\frac{n}{2}-D\left(A\right)$
with $A$ such that $|A|=\left\lceil a2^{n}\right\rceil $.  In the
bottom figure, the vertical axis corresponds to the degree-$1$ Fourier
weight $\mathbf{W}_{1}$ of a Boolean function $f$ such that $\left|f^{-1}(1)\right|=\left\lceil a2^{n}\right\rceil $.
The quantities $\frac{n}{2}-D\left(A\right)$ and $\mathbf{W}_{1}$
are related via \eqref{eq:-27}. The circles correspond to $\frac{1}{2}\log_{2}\frac{1}{a}$
with $a=2^{-i},i=1,2,3,...$, which are the gaps $\frac{n}{2}-D\left(A\right)$
for subcube codes $A=\{1\}^{i}\times\left\{ -1,1\right\} ^{n-i}$.
The ``Hamming Ball'' curve corresponds to the average distances
of the Hamming balls $\left\{ \mathbf{x}:d_{\mathrm{H}}\left(\mathbf{x},\mathbf{0}\right)\le i\right\} $
for $i\in[1:n]$. The average distances of the Hamming balls are characterized
in \cite[Proposition 5.25]{O'Donnell14analysisof}. Here we only plot
bounds for $0\le a\le\frac{1}{2}$, since by using the relationship
given in \eqref{eq:-11-1-1}, the bounds corresponding to $\frac{1}{2}<a\le1$
are implied by the bounds corresponding to $0\le a\le\frac{1}{2}$.}
\end{figure}

Fu-Wei-Yeung's bound in \eqref{eq:-Fu}, the improved linear programming
bound in \eqref{eq:-56}, Chang's bound in \eqref{eq:-31}, and the
hypercontractivity bound in \eqref{eq:-39} are plotted in Fig. \ref{fig:simulation-1}.
Our improved linear programming bound is tighter than Fu-Wei-Yeung's
bound for $a<1/4$. It is tighter than Chang's bound (resp. the hypercontractivity
bound) when $a$ is larger than a value of approximately $0.08$ (resp.
a value of approximately $0.09$). The average distances of Hamming
subcubes are smaller than those of Hamming balls when $a$ is large,
and larger than those of Hamming balls when $a$ is small. For $a=1/2$
or $1/4$, Hamming subcubes attain the minimum average distance. However,
if $a$ tends to zero, Hamming balls asymptotically attain the minimum
average distance among sets of volume $\left\lceil a2^{n}\right\rceil $
\cite[Remark 5.28]{O'Donnell14analysisof}. Our linear programming
bound is tighter than existing bounds for $a=1/8$. However, for this
case, there is still a gap between our lower bound and the average
distance of Hamming subcubes. The latter is the best known upper bound
on the minimum average distance for this case. Hence at present, it
is still unclear whether Hamming subcubes are optimal for $a=1/8$.

\section{Bounds on Fourier Weights}

By using the relationship \eqref{eq:-8} between $\mathbf{W}_{1}$
and $D\left(A\right)$ with $f=2\cdot1_{A}-1$, Theorem \ref{thm:ImprovedLPB}
implies the following bound on $\mathbf{W}_{1}$. In the following,
we denote $M=|f^{-1}(1)|$ and $a=\frac{M}{2^{n}}$ for a Boolean
function.
\begin{cor}[Improved LP Bound]
\label{thm:ImprovedLPB-2} For $a=\frac{M}{2^{n}}\le\frac{1}{2}$,
 the degree-$1$ Fourier weight of a Boolean function $f$ such that
$\left|f^{-1}(1)\right|=M$ satisfies 
\begin{equation}
\mathbf{W}_{1}\leq8a^{2}\varphi(a).\label{eq:-56-2}
\end{equation}
\end{cor}
\begin{rem}
The Hypercontractivity bound in Proposition \ref{prop:RHB-1} implies
that 
\begin{equation}
\mathbf{W}_{1}\leq8a^{2}\psi\left(a\right).\label{eq:-61}
\end{equation}
\end{rem}

We next upper bound $\mathbf{W}_{m}$ for $m\ge2$ by using a linear
programming approach similar to Fu-Wei-Yeung's approach. Since $\sum_{i=0}^{n}Q^{\left(A\right)}(i)=\frac{1}{a}$,
we have
\[
Q^{\left(A\right)}(m)=\frac{1}{a}-1-\sum_{i\neq0,m}Q^{\left(A\right)}(i)
\]
Hence 
\begin{align}
\mathbf{W}_{m} & =4a^{2}\left(\frac{1}{a}-1-\sum_{i\neq0,m}Q^{\left(A\right)}(i)\right).\label{eq:-83}
\end{align}
Now we consider the following related optimization problem.
\begin{problem}
\label{prob:Primal-Problem:-1-1}Primal Problem:
\begin{equation}
\Phi_{m}(n;a):=\min_{u_{1},...,u_{m-1},u_{m+1},..,u_{n}}\sum_{i\neq0,m}u_{i}\label{eq:-84}
\end{equation}
subject to the inequalities 
\begin{align*}
 & u_{i}\ge0,\;i\in\left[1:m-1\right]\cup\left[m+1:n\right];\\
 & \sum_{i=1,i\neq m}^{n}\left[K_{k}\left(m\right)-K_{k}\left(i\right)\right]u_{i}\le K_{k}\left(0\right)+\left(\frac{1}{a}-1\right)K_{k}\left(m\right),\;k\in\left[1:n\right].
\end{align*}
\end{problem}
The dual of Problem \ref{prob:Primal-Problem:-1-1} is given as follows.

\begin{problem}
\label{prob:Dual-Problem:-1-1}Dual Problem:
\begin{equation}
\overline{\Phi}_{m}(n;a):=\max_{x_{1},x_{2},...,x_{n}}-\sum_{k=1}^{n}\left[K_{k}\left(0\right)+\left(\frac{1}{a}-1\right)K_{k}\left(m\right)\right]x_{k}\label{eq:-85}
\end{equation}
subject to 
\begin{align*}
 & x_{k}\ge0,\;k\in\left[1:n\right];\\
 & \sum_{k=1}^{n}\left[K_{k}\left(m\right)-K_{k}\left(i\right)\right]x_{k}\ge-1,\;i\in\left[1:m-1\right]\cup\left[m+1:n\right].
\end{align*}
\end{problem}

By strong duality of linear programming (and feasibility and boundedness
of the primal problem),

\[
\sum_{i\neq0,m}Q^{\left(A\right)}(i)\geq\Phi(n;a)=\overline{\Phi}(n;a).
\]
Then we prove the following bounds on Fourier weights. Since the proof
of Theorem \ref{thm:ImprovedLPB-1} is similar to that of Theorem
\ref{thm:ImprovedLPB}, it is omitted for the sake of brevity.
\begin{thm}[Bounds on Degree-$m$ Fourier Weight]
\label{thm:ImprovedLPB-1} For $a=\frac{M}{2^{n}}\le\frac{1}{2}$,
let $\mathbf{W}_{m}$ be the degree-$m$ Fourier weight of a Boolean
function $f$ such that $\left|f^{-1}(1)\right|=M$. For even $m\ge2$,
\begin{align}
\mathbf{W}_{m} & \leq4a(1-a).\label{eq:-57}
\end{align}
For odd $m\ge3$, 
\begin{align}
\mathbf{W}_{m} & \leq2a.\label{eq:-54}
\end{align}
\end{thm}
In proving Theorem \ref{thm:ImprovedLPB-1}, we use the feasible solutions
$\left(0,0,...,0,0\right)$ and $\left(0,0,...,0,\frac{1}{2}\right)$
to Problem \ref{prob:Dual-Problem:-1-1} to prove \eqref{eq:-57}
and \eqref{eq:-54} respectively. In the following, we show that these
two solutions are asymptotically optimal. That is, by using the linear
programming approach in \eqref{eq:-84} and \eqref{eq:-85}, it is
not possible to obtain better asymptotic bound as $n\to\infty$ .
The proof of Theorem \ref{thm:Optimality-1} is similar to that of
Theorem \ref{thm:Optimality}, and hence is also omitted here.
\begin{thm}[Asymptotic Optimality of the Bounds in \eqref{eq:-54} and \eqref{eq:-57}]
\label{thm:Optimality-1} For fixed $a\le\frac{1}{2}$ and $m\in[1:n]$,
and for $M=\left\lceil a2^{n}\right\rceil $, 
\begin{equation}
\lim_{n\to\infty}\overline{\Phi}_{m}(n;a)=\begin{cases}
0 & \textrm{even }m\ge2\\
\frac{1}{2a}-1 & \textrm{odd }m\ge3
\end{cases}.\label{eq:-72-1-1}
\end{equation}
\end{thm}
Theorem \ref{thm:Optimality-1} implies that for the asymptotic case
as $n\to\infty$, the bounds in Theorem \ref{thm:ImprovedLPB-1} are
the best possible that can be obtained via the linear programming
approach (i.e., Problems \ref{prob:Primal-Problem:-1-1} and \ref{prob:Dual-Problem:-1-1}).
It is somewhat interesting to note that for $\mathbf{W}_{1}$, to
achieve asymptotic optimality, the feasible solution has to be $2$-sparse
(for $a\le1/4$). However, for $\mathbf{W}_{m}$ in which $m\ge2$,
a $1$-sparse or $0$-sparse dual feasible solution suffices for achieving
asymptotic optimality.

\section{Application to Noise Stability}

In this section, we apply our bounds on the degree-$1$ Fourier weight
to bound the noise stability of Boolean functions. Let
\begin{align}
 & \begin{array}{c}
\qquad\qquad\qquad-1\qquad1\\
P_{XY}=\begin{array}{c}
-1\\
1
\end{array}\left[\begin{array}{cc}
\frac{1+\rho}{4} & \frac{1-\rho}{4}\\
\frac{1-\rho}{4} & \frac{1+\rho}{4}
\end{array}\right]
\end{array}\label{eq:-34}
\end{align}
be a joint distribution on $\{-1,1\}$ with correlation coefficient
$\rho\in[-1,1]$. Let $(\mathbf{X},\mathbf{Y})\sim P_{XY}^{n}$ be
$n$ i.i.d. copies of $\left(X,Y\right)\sim P_{XY}$.  

\begin{defn}
\label{def:-For-two-1}  For $f:\{-1,1\}^{n}\to\{-1,1\}$ and $\rho\in[-1,1]$,
the \emph{noise stability} of $f$ at $\rho$ is 
\[
\mathbf{Stab}_{\rho}[f]:=\mathbb{E}[f(\mathbf{X})f(\mathbf{Y})],
\]
where the expectation is taken over random vectors $(\mathbf{X},\mathbf{Y})\sim P_{XY}^{n}$
with $P_{XY}$ defined in \eqref{eq:-34}. 
\end{defn}

For a Boolean function $f:\{-1,1\}^{n}\to\{-1,1\}$ such that $\mathbb{P}\left(f(\mathbf{X})=1\right)=a$,
define 
\[
q:=\mathbb{P}\left(f(\mathbf{X})=f(\mathbf{Y})=1\right).
\]
Then 
\begin{align*}
\mathbf{Stab}_{\rho}[f] & =2\mathbb{P}\left(f(\mathbf{X})=f(\mathbf{Y})\right)\text{\textendash}1\\
 & =2\left(1+2q-2a\right)\text{\textendash}1.
\end{align*}
Given $a$, maximizing $\mathbf{Stab}_{\rho}[f]$ is equivalent to
maximizing $q$. On the other hand, by \cite[Equation (11)]{yu2019bounds},
$q=a^{2}\Pi_{\rho}\left(A\right)$, where $\Pi_{\rho}\left(A\right)$
is the distance enumerator of $A$ defined in \eqref{eq:-51}. By
\eqref{eq:-52}, 
\begin{align}
q & =a^{2}\sum_{k\geq0}\rho^{k}\cdot Q^{\left(A\right)}(k)\\
 & =a^{2}+\frac{1}{4}\sum_{k\geq1}\rho^{k}\cdot\mathbf{W}_{k}.\label{eq:-60}
\end{align}
Since $\sum_{k\geq0}\mathbf{W}_{k}=1$, Equation \eqref{eq:-60} leads
to the following inequalities:
\begin{equation}
a^{2}-\frac{1}{4}\rho\cdot\mathbf{W}_{1}-\frac{1}{4}\rho^{2}\cdot\left(1-\left(2a-1\right)^{2}-\mathbf{W}_{1}\right)\leq q\leq a^{2}+\frac{1}{4}\rho\cdot\mathbf{W}_{1}+\frac{1}{4}\rho^{2}\cdot\left(1-\left(2a-1\right)^{2}-\mathbf{W}_{1}\right)\label{eq:-7-1-1}
\end{equation}
i.e.,
\begin{equation}
a^{2}-\frac{1}{4}\rho\cdot\mathbf{W}_{1}-\frac{1}{4}\rho^{2}\cdot\left(4a\left(1-a\right)-\mathbf{W}_{1}\right)\leq q\leq a^{2}+\frac{1}{4}\rho\cdot\mathbf{W}_{1}+\frac{1}{4}\rho^{2}\cdot\left(4a\left(1-a\right)-\mathbf{W}_{1}\right).\label{eq:-7-1-1-2}
\end{equation}
Define 
\[
\eta\left(a\right):=\min\left\{ \varphi(a),\psi\left(a\right)\right\} 
\]
then by the bounds on $\mathbf{W}_{1}$ given in \eqref{eq:-56-2}
and \eqref{eq:-61},
\begin{equation}
\mathbf{W}_{1}\leq8a^{2}\eta\left(a\right).\label{eq:-56-2-1}
\end{equation}
Therefore,
\begin{equation}
a^{2}-2a^{2}\eta\left(a\right)\rho-\left(a\left(1-a\right)-2a^{2}\eta\left(a\right)\right)\rho^{2}\leq q\leq a^{2}+2a^{2}\eta\left(a\right)\rho+\left(a\left(1-a\right)-2a^{2}\eta\left(a\right)\right)\rho^{2}.\label{eq:-7-1-1-1}
\end{equation}
Hence we have the following bounds on $q$.
\begin{prop}
\label{prop:symmetric-1} 
\begin{equation}
\theta^{-}(a)\leq q\leq\theta^{+}(a),\label{eq:-82}
\end{equation}
where 
\[
\theta^{+}(a):=\min\left\{ a,a^{2}+2a^{2}\eta\left(a\right)\rho+\left(a\left(1-a\right)-2a^{2}\eta\left(a\right)\right)\rho^{2}\right\} 
\]
and 
\[
\theta^{-}(a):=\max\left\{ 0,a^{2}-2a^{2}\eta\left(a\right)\rho-\left(a\left(1-a\right)-2a^{2}\eta\left(a\right)\right)\rho^{2}\right\} .
\]
\end{prop}
If we replace $\eta\left(a\right)$ with $\frac{1}{4a}$, then Proposition
\ref{prop:symmetric-1} reduces to \cite[Corollary 1]{yu2019bounds}.
It is easy to verify that $\varphi(a)\leq\frac{1}{4a}$. Hence $\eta\left(a\right)\leq\frac{1}{4a}$,
which implies that the lower and upper bounds given in Proposition
\ref{prop:symmetric-1} are tighter than the corresponding bounds
given in the present authors' previous work \cite[Corollary 1]{yu2019bounds}.

\appendices{}

\section{\label{sec:Proof-of-Proposition}Proof of Proposition \ref{prop:monotonicity}}

Let $A^{*}\subseteq\{-1,1\}^{n}$ be an $(n,M)$-code that attains
 $\min_{A\subseteq\{-1,1\}^{n}:|A|=M}D\left(A\right)$. Now we construct
a new $(n+k,2^{k}M)$-code as follows:
\[
B=A^{*}\times\{-1,1\}^{k}.
\]
Obviously, $B\subseteq\{-1,1\}^{n+k}$ and $|B|=2^{k}M$. Next, we
prove that $D\left(B\right)=D\left(A^{*}\right)$.

For any $\mathbf{x}\in B$, we can write $\mathbf{x}=\left(\mathbf{x}_{1},\mathbf{x}_{2}\right)$
where $\mathbf{x}_{1}\in A^{*}$ and $\mathbf{x}_{2}\in\{-1,1\}^{k}$.
Then we have
\begin{equation}
d_{\mathrm{H}}\left(\mathbf{x},\mathbf{x}'\right)=d_{\mathrm{H}}\left(\mathbf{x}_{1},\mathbf{x}_{1}'\right)+d_{\mathrm{H}}\left(\mathbf{x}_{2},\mathbf{x}_{2}'\right).\label{eq:-68}
\end{equation}
Using \eqref{eq:-68} we obtain that 
\begin{align}
D\left(B\right) & =\frac{1}{|B|^{2}}\sum_{\left(\mathbf{x},\mathbf{x}'\right)\in B^{2}}d_{\mathrm{H}}\left(\mathbf{x},\mathbf{x}'\right)\nonumber \\
 & =\frac{1}{\left(2^{k}M\right)^{2}}\sum_{\left(\mathbf{x},\mathbf{x}'\right)\in B^{2}}\left[d_{\mathrm{H}}\left(\mathbf{x}_{1},\mathbf{x}_{1}'\right)+d_{\mathrm{H}}\left(\mathbf{x}_{2},\mathbf{x}_{2}'\right)\right]\nonumber \\
 & =\frac{1}{\left(2^{k}M\right)^{2}}\sum_{\left(\mathbf{x},\mathbf{x}'\right)\in B^{2}}d_{\mathrm{H}}\left(\mathbf{x}_{1},\mathbf{x}_{1}'\right)+\frac{1}{\left(2^{k}M\right)^{2}}\sum_{\left(\mathbf{x},\mathbf{x}'\right)\in B^{2}}d_{\mathrm{H}}\left(\mathbf{x}_{2},\mathbf{x}_{2}'\right)\nonumber \\
 & =\frac{2^{2k}}{\left(2^{k}M\right)^{2}}\sum_{\left(\mathbf{x}_{1},\mathbf{x}_{1}'\right)\in A^{*2}}d_{\mathrm{H}}\left(\mathbf{x}_{1},\mathbf{x}_{1}'\right)+\frac{M^{2}}{\left(2^{k}M\right)^{2}}\sum_{\left(\mathbf{x}_{2},\mathbf{x}_{2}'\right)\in\{-1,1\}^{2k}}d_{\mathrm{H}}\left(\mathbf{x}_{2},\mathbf{x}_{2}'\right)\nonumber \\
 & =D\left(A^{*}\right)+\frac{k}{2},\label{eq:-66}
\end{align}
where \eqref{eq:-66} follows since for a $k$-dimensional Hamming
cube, its average distance is $k/2$.

Hence 
\[
\min_{A\subseteq\{-1,1\}^{n+k}:|A|=2^{k}M}D\left(A\right)\leq D\left(B\right)=D\left(A^{*}\right)+\frac{k}{2}.
\]

\section{\label{sec:Proof-of-Theorem-ImprovedLPB}Proof of Theorem \ref{thm:ImprovedLPB}}

We first provide some fundamental properties of Krawtchouk polynomials
in Appendix \ref{subsec:Properties-of-Krawtchouk}, and then applied
them to prove Theorem \ref{thm:ImprovedLPB} in Appendix \ref{subsec:Proof-of-Theorem}.

\subsection{\label{subsec:Properties-of-Krawtchouk}Properties of Krawtchouk
Polynomials}

By definition, Krawtchouk polynomials satisfy 
\begin{align}
K_{k}\left(n-i\right) & =(-1)^{k}K_{k}\left(i\right)\label{eq:-78}\\
K_{n-k}\left(i\right) & =(-1)^{i}K_{k}\left(i\right);\label{eq:-79}
\end{align}
see \cite{macwilliams1977theory}. Furthermore, Krawtchouk polynomials
also have the following recurrence property.
\begin{lem}
\label{lem:}For $x\in[0,n-1]$,
\[
K_{k}^{(n)}\left(x+1\right)=K_{k}^{(n)}\left(x\right)-2K_{k-1}^{(n-1)}\left(x\right).
\]
\end{lem}
\begin{IEEEproof}
By \cite[Theorem 15]{macwilliams1977theory}, we have the following
alternative expression for Krawtchouk polynomials:
\begin{equation}
K_{k}^{(n)}\left(x\right)=\sum_{j=0}^{k}(-2)^{j}{x \choose j}{n-j \choose k-j}.\label{eq:-62}
\end{equation}
By using the alternative expression in \eqref{eq:-62}, we obtain
that 
\begin{align*}
K_{k}^{(n)}\left(x\right)-K_{k}^{(n)}\left(x+1\right) & =\sum_{j=0}^{k}(-2)^{j}{x \choose j}{n-j \choose k-j}-\sum_{j=0}^{k}(-2)^{j}{x+1 \choose j}{n-j \choose k-j}\\
 & =\sum_{j=0}^{k}(-2)^{j}{x \choose j}{n-j \choose k-j}\left(1-\frac{x+1}{x+1-j}\right)\\
 & =-\sum_{j=0}^{k}(-2)^{j}{x \choose j}{n-j \choose k-j}\frac{j}{x+1-j}\\
 & =-\sum_{j=0}^{k}(-2)^{j}{x \choose j-1}{n-j \choose k-j}\\
 & =2\sum_{j=0}^{k}(-2)^{j-1}{x \choose j-1}{\left(n-1\right)-\left(j-1\right) \choose \left(k-1\right)-\left(j-1\right)}\\
 & =2\sum_{j'=0}^{k-1}(-2)^{j'}{x \choose j'}{\left(n-1\right)-j' \choose \left(k-1\right)-j'}\\
 & =2K_{k-1}^{(n-1)}\left(x\right)
\end{align*}
\end{IEEEproof}
We consider the function $\mathbb{R}\ni x\mapsto K_{k}^{(n)}\left(x\right)\in\mathbb{R}$,
which has $k$ distinct real roots \cite{macwilliams1977theory}.
We denote the real roots respectively as $x_{1}^{(n,k)}<x_{2}^{(n,k)}<...<x_{k}^{(n,k)}$.
These roots lie in the interval $\left[\frac{n}{2}-\sqrt{k\left(n-k\right)},\frac{n}{2}+\sqrt{k\left(n-k\right)}\right]$
and are symmetric with respect to the point $x=\frac{n}{2}$ (i.e.,
$x_{i}^{(n,k)}+x_{k+1-i}^{(n,k)}=n$ for $1\le i\le k$) \cite[Section 2.1.2]{kirshner2019moment}.
Moreover, $K_{k}^{(n)}\left(x\right)\geq0$ for $0\leq x\leq x_{1}^{(n,k)}$.
Hence we have the following properties.
\begin{lem}
\label{lem:The-following-hold:}The following hold:
\begin{enumerate}
\item (Monotonicity) For $0\leq x\leq x_{1}^{(n-1,k-1)}$, we have 
\begin{equation}
K_{k}^{(n)}\left(x\right)\geq K_{k}^{(n)}\left(x+1\right).\label{eq:-63}
\end{equation}
\item (Bound on Magnitude) For all $x\in[0:n]$,
\begin{equation}
\left|K_{k}^{(n)}\left(x\right)\right|<2^{\frac{n}{2}\left(1+H\left(\frac{k}{n}\right)-H\left(\frac{x}{n}\right)+\frac{1}{n}\log_{2}(n+1)\right)}.\label{eq:-86}
\end{equation}
\end{enumerate}
\end{lem}
\begin{IEEEproof}
Statement 1 follows by Lemma \ref{lem:} and the fact $K_{k-1}^{(n-1)}\left(x\right)\geq0$
for $0\leq x\leq x_{1}^{(n-1,k-1)}$. Now we prove Statement 2. By
\cite[Equation (16)]{krasikov2001nonnegative}, for $x\in[0:n]$,
\[
K_{k}^{(n)}\left(x\right)^{2}<2^{n}{n \choose k}{n \choose x}^{-1}.
\]
By \cite[Lemma 2.3 and Problem 1 on p.39]{Csiszar},
\[
\frac{1}{n+1}2^{nH\left(\frac{k}{n}\right)}\leq{n \choose k}\leq2^{nH\left(\frac{k}{n}\right)}.
\]
Hence for all $x\in[0:n]$,
\[
\left|K_{k}^{(n)}\left(x\right)\right|<2^{\frac{n}{2}\left(1+H\left(\frac{k}{n}\right)-H\left(\frac{x}{n}\right)+\frac{1}{n}\log_{2}(n+1)\right)}.
\]
\end{IEEEproof}
\begin{lem}
\label{lem:Ki_is_larger}Let $i\in\mathbb{Z}_{\ge0}$ and $0<\beta<1/2$.
Let $k=\left\lfloor \beta n\right\rfloor $. Then given $i$ and $\beta$,
there exists an $N_{i,\beta}\in\mathbb{Z}_{\ge0}$ such that for all
$n\ge N_{i,\beta}$, 
\begin{equation}
K_{k}^{(n)}\left(i\right)\geq\left|K_{k}^{(n)}\left(x\right)\right|,\quad\forall x\in\left[i,n-i\right].\label{eq:-64}
\end{equation}
\end{lem}

\begin{IEEEproof}
By \eqref{eq:-78}, the function $x\mapsto\left|K_{k}^{(n)}\left(x\right)\right|$
is symmetric with respect to the line $x=\frac{n}{2}$. Hence to prove
that \eqref{eq:-64} holds for $i\leq x\leq n-i$, it suffices to
prove that it holds for $i\leq x\leq\frac{n}{2}$. Next we prove this.

For $i\leq x\leq x_{1}^{(n-1,k-1)}$, \eqref{eq:-64} follows by
\eqref{eq:-63}. Now we consider $x_{1}^{(n-1,k-1)}\le x\leq\frac{n}{2}$.

By \eqref{eq:-86}, we have that 
\begin{align}
\frac{K_{k}^{(n)}\left(i\right)}{\left|K_{k}^{(n)}\left(x\right)\right|} & >\frac{K_{k}^{(n)}\left(i\right)}{2^{\frac{n}{2}\left(1+H\left(\frac{k}{n}\right)-H\left(\frac{x}{n}\right)+\frac{1}{n}\log_{2}(n+1)\right)}}\nonumber \\
 & =\frac{{n \choose k}}{2^{\frac{n}{2}\left(1+H\left(\frac{k}{n}\right)-H\left(\frac{x}{n}\right)+\frac{1}{n}\log_{2}(n+1)\right)}}\frac{K_{k}^{(n)}\left(i\right)}{{n \choose k}}\nonumber \\
 & \geq\frac{(n+1)^{-1}2^{nH\left(\frac{k}{n}\right)}}{2^{\frac{n}{2}\left(1+H\left(\frac{k}{n}\right)-H\left(\frac{x}{n}\right)+\frac{1}{n}\log_{2}(n+1)\right)}}\frac{K_{k}^{(n)}\left(i\right)}{{n \choose k}}\\
 & =2^{\frac{n}{2}\left(-1+H\left(\frac{k}{n}\right)+H\left(\frac{x}{n}\right)-\frac{3}{n}\log_{2}(n+1)\right)}\frac{K_{k}^{(n)}\left(i\right)}{{n \choose k}}.\label{eq:-47}
\end{align}
Since by \eqref{eq:-62}, for $k\ge i$, $K_{k}^{(n)}\left(i\right)=\sum_{j=0}^{i}(-2)^{j}{i \choose j}{n-j \choose k-j}$,
we have 
\begin{align}
\frac{K_{k}^{(n)}\left(i\right)}{{n \choose k}} & =\frac{\sum_{j=0}^{i}(-2)^{j}{i \choose j}{n-j \choose k-j}}{{n \choose k}}\label{eq:-62-1}\\
 & =\sum_{j=0}^{i}(-2)^{j}\frac{k\left(k-1\right)\cdots\left(k-j+1\right)}{n\left(n-1\right)\cdots\left(n-j+1\right)}{i \choose j}\\
 & \to\sum_{j=0}^{i}(-2)^{j}\beta^{j}{i \choose j}\\
 & =\left(1-2\beta\right)^{i}>0.\label{eq:-65}
\end{align}
On the other hand, for all $x_{1}^{(n-1,k-1)}\le x\leq\frac{n}{2}$,
\begin{align}
H\left(\frac{k}{n}\right)+H\left(\frac{x}{n}\right) & \geq H\left(\frac{k}{n}\right)+H\left(\frac{1}{2}-\sqrt{\frac{k-1}{n-1}\left(1-\frac{k-1}{n-1}\right)}\right)\nonumber \\
 & \to H\left(\beta\right)+H\left(\frac{1}{2}-\sqrt{\beta\left(1-\beta\right)}\right)\nonumber \\
 & >1,\label{eq:-45}
\end{align}
where \eqref{eq:-45} follows from \cite[Section 2.1.2]{kirshner2019moment}
and the hypothesis $0<\beta<1/2$ . Combining \eqref{eq:-65} and
\eqref{eq:-45} yields that \eqref{eq:-47} is exponentially large.
This means that for sufficiently large $n$, $K_{k}^{(n)}\left(i\right)\geq\left|K_{k}^{(n)}\left(x\right)\right|$
holds for all $x_{1}^{(n-1,k-1)}\le x\leq\frac{n}{2}$.
\end{IEEEproof}

\subsection{\label{subsec:Proof-of-Theorem}Proof of Theorem \ref{thm:ImprovedLPB}}

By Proposition \ref{prop:monotonicity}, 
\begin{align*}
D\left(A\right)-\frac{n}{2} & \geq\liminf_{n\to\infty}D\left(A\right)-\frac{n}{2}\\
 & \geq\frac{1}{2}-\frac{1}{2a}+\frac{1}{2}\liminf_{n\to\infty}\overline{\Lambda}(n;a).
\end{align*}
Hence to prove Theorem \ref{thm:ImprovedLPB}, it suffices to prove
\begin{align}
\liminf_{n\to\infty}\overline{\Lambda}(n;a) & \geq\begin{cases}
\frac{\left(1-\sqrt{a}\right)^{2}}{a} & 0\le a<1/4\\
\frac{1}{2a}-1 & 1/4\le a\le1/2
\end{cases}.\label{eq:-73}
\end{align}
Next we prove this. 

Let $\beta\in\left(\frac{1}{2},1\right)$ be a constant. Let $k=2\left\lfloor \frac{\beta n}{2}\right\rfloor $.
Obviously, $k$ is an even integer and $k/n\to\beta$ as $n\to\infty$.
Then we consider the vector $\mathbf{x}^{*}:=\left(0,...,0,x_{k}^{*},x_{k+1}^{*},0,...,0\right)$
with the $k$-th and $\left(k+1\right)$-th components $\left(x_{k}^{*},x_{k+1}^{*}\right)$
satisfying 
\begin{align}
 & \left[K_{k}\left(2\right)-K_{k}\left(1\right)\right]x_{k}^{*}+\left[K_{k+1}\left(2\right)-K_{k+1}\left(1\right)\right]x_{k+1}^{*}=1\label{eq:-80}\\
 & \left[K_{k}\left(n\right)-K_{k}\left(1\right)\right]x_{k}^{*}+\left[K_{k+1}\left(n\right)-K_{k+1}\left(1\right)\right]x_{k+1}^{*}=1.\label{eq:-81}
\end{align}
Solving the equations \eqref{eq:-80} and \eqref{eq:-81}, we obtain
that 
\begin{align*}
x_{k}^{*} & =\frac{K_{k+1}\left(2\right)-K_{k+1}\left(n\right)}{\left[K_{k+1}\left(2\right)-K_{k+1}\left(1\right)\right]\left[K_{k}\left(n\right)-K_{k}\left(1\right)\right]-\left[K_{k}\left(2\right)-K_{k}\left(1\right)\right]\left[K_{k+1}\left(n\right)-K_{k+1}\left(1\right)\right]}\\
x_{k+1}^{*} & =\frac{K_{k}\left(n\right)-K_{k}\left(2\right)}{\left[K_{k+1}\left(2\right)-K_{k+1}\left(1\right)\right]\left[K_{k}\left(n\right)-K_{k}\left(1\right)\right]-\left[K_{k}\left(2\right)-K_{k}\left(1\right)\right]\left[K_{k+1}\left(n\right)-K_{k+1}\left(1\right)\right]}.
\end{align*}
Observe that 
\begin{align*}
K_{k}\left(n\right) & =(-1)^{k}K_{k}\left(0\right)=(-1)^{k}{n \choose k},\\
K_{k}\left(1\right) & ={n \choose k}\left(1-\frac{2k}{n}\right),\\
K_{k}\left(2\right) & =\sum_{j=0}^{k}(-1)^{j}{2 \choose j}{n-2 \choose k-j}\\
 & ={n-2 \choose k}-2{n-2 \choose k-1}+{n-2 \choose k-2}\\
 & ={n \choose k}\left(\frac{\left(n-k\right)\left(n-k-1\right)}{n\left(n-1\right)}-2\frac{\left(n-k\right)k}{n\left(n-1\right)}+\frac{k\left(k-1\right)}{n\left(n-1\right)}\right)\\
 & ={n \choose k}\frac{\left(n-2k\right)^{2}-n}{n\left(n-1\right)}.
\end{align*}
Therefore,
\begin{align*}
x_{k}^{*} & =\frac{1+2\left(\frac{k}{n}+\frac{1}{n}\right)^{2}-2\left(\frac{k}{n}+\frac{1}{n}\right)-\frac{1}{n}}{{n \choose k}2\frac{k}{n}\left(2\frac{k}{n}-1+\frac{1}{n}\left(2\frac{k}{n}+2\frac{1}{n}-1\right)\right)}\\
x_{k+1}^{*} & =\frac{1-\frac{k}{n}}{{n \choose k+1}\left(2\frac{k}{n}-1+\frac{1}{n}\left(2\frac{k}{n}+2\frac{1}{n}-1\right)\right)}
\end{align*}
as stated in \eqref{eq:oursolution} and \eqref{eq:oursolution-1}.
Letting $n\to\infty$, we have 
\begin{align*}
{n \choose k}x_{k}^{*} & \to\frac{1+2\beta^{2}-2\beta}{2\beta\left(2\beta-1\right)}=\frac{1}{2\beta\left(2\beta-1\right)}-\frac{1-\beta}{2\beta-1}\\
{n \choose k+1}x_{k+1}^{*} & \to\frac{1-\beta}{2\beta-1}.
\end{align*}
We next prove that when $\beta>1/2$ and $n$ is sufficiently large,
$x^{*}=\left(0,...,0,x_{k}^{*},x_{k+1}^{*},0,...,0\right)$ is a feasible
solution to Problem \ref{prob:Dual-Problem:}.

Observe that for sufficiently large $n$, $x_{k}^{*},x_{k+1}^{*}\geq0$.
Hence we only need to show that for sufficiently large $n$ and for
all $i=\left[2:n\right]$, the other inequality constraint in Problem
\ref{prob:Dual-Problem:} is satisfied, i.e.,
\begin{equation}
\varphi(i):=\left[K_{k}\left(i\right)-K_{k}\left(1\right)\right]x_{k}^{*}+\left[K_{k+1}\left(i\right)-K_{k+1}\left(1\right)\right]x_{k+1}^{*}\leq1.\label{eq:-29}
\end{equation}
By the choice of $\mathbf{x}^{*}$, we have $\varphi(2)=\varphi(n)=1$.
By \eqref{eq:-79}, we have 
\begin{align*}
\varphi(i) & =\left[K_{k}\left(i\right)-K_{k}\left(1\right)\right]x_{k}^{*}+\left[K_{k+1}\left(i\right)-K_{k+1}\left(1\right)\right]x_{k+1}^{*}\\
 & =\left[(-1)^{i}K_{n-k}\left(i\right)+K_{n-k}\left(1\right)\right]x_{k}^{*}+\left[(-1)^{i}K_{n-k-1}\left(i\right)+K_{n-k-1}\left(1\right)\right]x_{k+1}^{*}.
\end{align*}
By Lemma \ref{lem:Ki_is_larger}, for $i\in[2:n-2]$,
\begin{align*}
\varphi(i) & \leq\left[\left|K_{n-k}\left(i\right)\right|+K_{n-k}\left(1\right)\right]x_{k}^{*}+\left[\left|K_{n-k-1}\left(i\right)\right|+K_{n-k-1}\left(1\right)\right]x_{k+1}^{*}\\
 & \leq\left[K_{n-k}\left(2\right)+K_{n-k}\left(1\right)\right]x_{k}^{*}+\left[K_{n-k-1}\left(2\right)+K_{n-k-1}\left(1\right)\right]x_{k+1}^{*}\\
 & =\varphi(2)=1.
\end{align*}
 Hence it remains to verify that $\varphi(n-1)\leq1$ for sufficiently
large $n$. Consider, 
\begin{align*}
\varphi(n-1) & =\left[K_{k}\left(n-1\right)-K_{k}\left(1\right)\right]x_{k}^{*}+\left[K_{k+1}\left(n-1\right)-K_{k+1}\left(1\right)\right]x_{k+1}^{*}\\
 & =2\left|K_{k+1}\left(1\right)\right|x_{k+1}^{*}\\
 & =2{n \choose k+1}\left|\left(1-\frac{2\left(k+1\right)}{n}\right)\right|x_{k+1}^{*}\\
 & \rightarrow2\left(1-\beta\right)<1.
\end{align*}

Until now, we have shown that when $\beta>1/2$ and $n$ is sufficiently
large, $x^{*}:=\left(0,...,0,x_{k}^{*},x_{k+1}^{*},0,...,0\right)$
is a feasible solution to Problem \ref{prob:Dual-Problem:}. This
immediately yields the following bounds on Problem \ref{prob:Primal-Problem:}.
For sufficiently large $n$,
\[
\overline{\Lambda}(n;a)\geq-{n \choose k}\left[1+\left(1-\frac{2k}{n}\right)\left(\frac{1}{a}-1\right)\right]x_{k}^{*}-{n \choose k+1}\left[1+\left(1-\frac{2\left(k+1\right)}{n}\right)\left(\frac{1}{a}-1\right)\right]x_{k+1}^{*}.
\]
Taking limits as $n\to\infty$, we obtain that 
\begin{align}
\liminf_{n\to\infty}\overline{\Lambda}(n;a) & \geq-\left[1+\left(1-2\beta\right)\left(\frac{1}{a}-1\right)\right]\frac{1+2\beta^{2}-2\beta+1-\left(2\beta-1\right)^{2}}{2\beta\left(2\beta-1\right)}\nonumber \\
 & =-\left[1+\left(1-2\beta\right)\left(\frac{1}{a}-1\right)\right]\frac{1}{2\beta\left(2\beta-1\right)}.\label{eq:-48}
\end{align}
Since $\beta\in\left(\frac{1}{2},1\right)$ is arbitrary, we can
maximize the bound \eqref{eq:-48} over all $\beta\in\left(\frac{1}{2},1\right)$.
This yields that 
\begin{align*}
\liminf_{n\to\infty}\overline{\Lambda}(n;a) & \geq\sup_{\beta\in\left(\frac{1}{2},1\right)}-\left[1+\left(1-2\beta\right)\left(\frac{1}{a}-1\right)\right]\frac{1}{2\beta\left(2\beta-1\right)}\\
 & =\begin{cases}
\frac{\left(1-\sqrt{a}\right)^{2}}{a} & 0\le a<1/4\\
\frac{1}{2a}-1 & 1/4\le a\le1/2
\end{cases}
\end{align*}
where the optimal value of $\beta$ is 
\[
\beta^{*}=\begin{cases}
\frac{1}{2\left(1-\sqrt{a}\right)} & 0\le a<1/4\\
1 & 1/4\le a\le1/2
\end{cases}.
\]
 The proof of \eqref{eq:-73} is complete.

\section{\label{sec:Proof-of-Theorem-Optimality}Proof of Proposition \ref{thm:Optimality-2}}

Consider Problem \ref{prob:Dual-Problem:}. Any feasible solution
to Problem \ref{prob:Dual-Problem:} satisfies 
\begin{align*}
 & x_{k}\ge0,\;k\in\left[1:n\right];\\
 & \sum_{k=1}^{n}\left[K_{k}\left(i\right)-K_{k}\left(1\right)\right]x_{k}\leq1,\;i\in\left[2:n\right].
\end{align*}
Taking $i=2$, we have
\begin{align}
 & \sum_{k=1}^{n}\left[K_{k}\left(2\right)-K_{k}\left(1\right)\right]x_{k}\leq1.\label{eq:-49}
\end{align}
Since 
\begin{align*}
K_{k}\left(1\right) & ={n \choose k}\left(1-\frac{2k}{n}\right)\\
K_{k}\left(2\right) & ={n \choose k}\frac{\left(n-2k\right)^{2}-n}{n\left(n-1\right)}
\end{align*}
we have that \eqref{eq:-49} is equivalent to 
\begin{equation}
\sum_{k=1}^{n}\left[\frac{2k\left(2k-n-1\right)}{n\left(n-1\right)}\right]y_{k}\leq1\label{eq:-69}
\end{equation}
where 
\[
y_{k}:={n \choose k}x_{k}\ge0,\;k\in\left[1:n\right].
\]
 The objective function of Problem \ref{prob:Dual-Problem:} satisfies
\begin{align}
 & -\sum_{k=1}^{n}\left[K_{k}\left(0\right)+K_{k}\left(1\right)\left(\frac{1}{a}-1\right)\right]x_{k}\nonumber \\
 & =\sum_{k=1}^{n}\left[-\frac{1}{a}+\frac{2k}{n}\left(\frac{1}{a}-1\right)\right]y_{k}\nonumber \\
 & \leq\sum_{k=\left\lceil \frac{n}{2\left(1-a\right)}\right\rceil }^{n}\left[-\frac{1}{a}+\frac{2k}{n}\left(\frac{1}{a}-1\right)\right]y_{k}.\label{eq:-70}
\end{align}
For $n\ge1/a-1$, we have $\frac{n}{2\left(1-a\right)}\ge\frac{n+1}{2}$.
Therefore, for $n\ge1/a-1$, the coefficients $\frac{2k\left(2k-n-1\right)}{n\left(n-1\right)}$
in \eqref{eq:-69} are nonnegative for all $k\in\left[\left\lceil \frac{n}{2\left(1-a\right)}\right\rceil :n\right]$.
By this property and continuing the upper bound in \eqref{eq:-70},
we obtain that for $n\ge1/a-1$,
\begin{align}
 & \sum_{k=\left\lceil \frac{n}{2\left(1-a\right)}\right\rceil }^{n}\left[-\frac{1}{a}+\frac{2k}{n}\left(\frac{1}{a}-1\right)\right]y_{k}\nonumber \\
 & \leq\left(\max_{k\in\left[\left\lceil \frac{n}{2\left(1-a\right)}\right\rceil :n\right]}\frac{-\frac{1}{a}+\frac{2k}{n}\left(\frac{1}{a}-1\right)}{\frac{2k\left(2k-n-1\right)}{n\left(n-1\right)}}\right)\sum_{k=\left\lceil \frac{n}{2\left(1-a\right)}\right\rceil }^{n}\left[\frac{2k\left(2k-n-1\right)}{n\left(n-1\right)}\right]y_{k}\label{eq:-50}\\
 & \leq\max_{k\in\left[\left\lceil \frac{n}{2\left(1-a\right)}\right\rceil :n\right]}\frac{-\frac{1}{a}+\frac{2k}{n}\left(\frac{1}{a}-1\right)}{\frac{2k\left(2k-n-1\right)}{n\left(n-1\right)}}\label{eq:-53}\\
 & \leq\max_{t\in\left[\left\lceil \frac{n}{2\left(1-a\right)}\right\rceil ,n\right]}\frac{-\frac{1}{a}+\frac{2t}{n}\left(\frac{1}{a}-1\right)}{\frac{2t\left(2t-n-1\right)}{n\left(n-1\right)}}\label{eq:-67}\\
 & =:\theta_{n}(a)
\end{align}
where \eqref{eq:-50} follows from the following inequality
\[
\sum_{i=1}^{m}a_{i}\le\left(\sum_{i=1}^{m}b_{i}\right)\max_{1\le i\le m}\left\{ \frac{a_{i}}{b_{i}}\right\} 
\]
for $b_{i}\ge0$ and real $a_{i},i\in\left[1:m\right]$; \eqref{eq:-53}
follows from \eqref{eq:-69}; and in \eqref{eq:-67}, the integer-valued
variable $k$ is relaxed to a real-valued variable $t$. 

Now we calculate the value of $\theta_{n}(a)$. By setting the derivative
of the objective function in \eqref{eq:-67} to be zero, we find that
for $n\ge1/a-1$, the objective function has a local minimum at 
\[
t_{1}=\frac{n-\sqrt{n\left(an+a-1\right)}}{2\left(1-a\right)}.
\]
and a local maximum at 
\[
t_{2}=\frac{n+\sqrt{n\left(an+a-1\right)}}{2\left(1-a\right)}.
\]
It is easy to verify that $t_{1}$ and $t_{2}$ satisfy the following
properties.
\begin{enumerate}
\item For $n\ge1/a-1$,
\[
t_{1}\leq\frac{n}{2\left(1-a\right)}.
\]
\item For $a\leq1/4$ and $n\ge1/a-1$,
\[
\frac{n}{2}\leq t_{2}\leq n.
\]
\item For $a>1/4$ and $n\ge\frac{1-a}{a-\left(1-2a\right)^{2}}$,
\[
t_{2}\geq n.
\]
\end{enumerate}
Based on the properties above, we know that the maximum in \eqref{eq:-67}
is attained at $t_{2}$ if $a\leq1/4$ and $n\ge1/a-1$, and at $n$
if $a>1/4$ and $n\ge\frac{1-a}{a-\left(1-2a\right)^{2}}$. That is,
for $a\leq1/4$ and $n\ge1/a-1$, 
\begin{align*}
\theta_{n}(a) & =\frac{-\frac{1}{a}+\frac{2t_{2}}{n}\left(\frac{1}{a}-1\right)}{\frac{2t_{2}\left(2t_{2}-n-1\right)}{n\left(n-1\right)}}\\
 & =\frac{\left(1-\frac{1}{n}\right)\frac{s}{a}}{\frac{1+s}{1-a}\left(\frac{1+s}{1-a}-1-\frac{1}{n}\right)}
\end{align*}
with $s:=\sqrt{a-\frac{1-a}{n}}$; and for $a>1/4$ and $n\ge\frac{1-a}{a-\left(1-2a\right)^{2}}$,
\begin{align*}
\theta_{n}(a) & =\frac{-\frac{1}{a}+2\left(\frac{1}{a}-1\right)}{2}\\
 & =\frac{1}{2a}-1.
\end{align*}
These yield \eqref{eq:-75} and \eqref{eq:-74} respectively.

\subsection*{Acknowledgements}

The authors are supported by  a Singapore Ministry of Education Tier
2 Grant (R-263-000-C83-112). \bibliographystyle{unsrt}
\bibliography{ref}

\begin{thebibliography}{10}

\bibitem{ahlswede1977contributions}
R.~Ahlswede and G.~O.~H. Katona.
\newblock Contributions to the geometry of {Hamming} spaces.
\newblock {\em Discrete Mathematics}, 17(1), 1977.

\bibitem{kundgen1998covering}
A.~K{\"u}ndgen.
\newblock Covering cliques with spanning bicliques.
\newblock {\em Journal of Graph Theory}, 27(4):223--227, 1998.

\bibitem{ahlswede1994asymptotic}
R.~Ahlswede and I.~Alth{\"o}fer.
\newblock The asymptotic behavior of diameters in the average.
\newblock {\em Journal of Combinatorial Theory, Series B}, 61(2):167--177,
  1994.

\bibitem{mounits2008lower}
B.~Mounits.
\newblock Lower bounds on the minimum average distance of binary codes.
\newblock {\em Discrete Mathematics}, 308(24):6241--6253, 2008.

\bibitem{althofer1992average}
I.~Alth{\"o}fer and T.~Sillke.
\newblock An ``average distance'' inequality for large subsets of the cube.
\newblock {\em Journal of Combinatorial Theory, Series B}, 56(2):296--301,
  1992.

\bibitem{shutao1998average}
S.~Xia and F.-W. Fu.
\newblock On the average {Hamming} distance for binary codes.
\newblock {\em Discrete Applied Mathematics}, 89(1-3):269--276, 1998.

\bibitem{fu1997expectation}
F.-W. Fu and S.-Y. Shen.
\newblock On the expectation and variance of {Hamming} distance between two iid
  random vectors.
\newblock {\em Acta Mathematicae Applicatae Sinica}, 13(3):243--250, 1997.

\bibitem{fu1999hamming}
F.-W. Fu, T.~Klove, and S.-Y. Shen.
\newblock On the {Hamming} distance between two iid random n-tuples over a
  finite set.
\newblock {\em IEEE Trans. Inf. Theory}, 45(2):803--807, 1999.

\bibitem{fu2001minimum}
F.-W. Fu, V.~K. Wei, and R.~W. Yeung.
\newblock On the minimum average distance of binary codes: {Linear} programming
  approach.
\newblock {\em Discrete Applied Mathematics}, 111(3):263--281, 2001.

\bibitem{mceliece1977new}
R.~McEliece, E.~Rodemich, H.~Rumsey, and L.~Welch.
\newblock New upper bounds on the rate of a code via the {Delsarte-MacWilliams}
  inequalities.
\newblock {\em IEEE Trans. Inf. Theory}, 23(2):157--166, 1977.

\bibitem{samorodnitsky2001optimum}
A.~Samorodnitsky.
\newblock On the optimum of {Delsarte's} linear program.
\newblock {\em Journal of Combinatorial Theory, Series A}, 96(2):261--287,
  2001.

\bibitem{yu2019bounds}
L.~Yu and V.~Y.~F. Tan.
\newblock Bounds on the average distance and distance enumerator with
  applications to non-interactive simulation.
\newblock {\em arXiv preprint arXiv:1904.03932}, 2019.

\bibitem{O'Donnell14analysisof}
R.~O'Donnell.
\newblock {\em Analysis of {Boolean} Functions}.
\newblock Cambridge University Press, 2014.

\bibitem{green2008boolean}
B.~Green and T.~Sanders.
\newblock Boolean functions with small spectral norm.
\newblock {\em Geometric and Functional Analysis}, 18(1):144--162, 2008.

\bibitem{chang2002polynomial}
M.-C. Chang.
\newblock A polynomial bound in {Freiman's} theorem.
\newblock {\em Duke mathematical journal}, 113(3):399--419, 2002.

\bibitem{friedgut2002boolean}
E.~Friedgut, G.~Kalai, and A.~Naor.
\newblock Boolean functions whose {Fourier} transform is concentrated on the
  first two levels.
\newblock {\em Advances in Applied Mathematics}, 29(3):427--437, 2002.

\bibitem{defant2019fourier}
A.~Defant, M.~Masty{\l}o, and A.~P{\'e}rez.
\newblock On the {Fourier} spectrum of functions on boolean cubes.
\newblock {\em Mathematische Annalen}, 374(1-2):653--680, 2019.

\bibitem{macwilliams1977theory}
F.~J. MacWilliams and N.~J.~A. Sloane.
\newblock {\em The Theory of Error-Correcting Codes}, volume~16.
\newblock Elsevier, 1977.

\bibitem{impagliazzo2014entropic}
R.~Impagliazzo, C.~Moore, and A.~Russell.
\newblock An entropic proof of {Chang's} inequality.
\newblock {\em SIAM Journal on Discrete Mathematics}, 28(1):173--176, 2014.

\bibitem{hambardzumyan2019chang}
L.~Hambardzumyan and Y.~Li.
\newblock Chang's lemma via {Pinsker's} inequality.
\newblock {\em Discrete Mathematics}, 2019.

\bibitem{kirshner2019moment}
N.~Kirshner and A.~Samorodnitsky.
\newblock A moment ratio bound for polynomials and some extremal properties of
  {Krawchouk} polynomials and {Hamming} spheres.
\newblock {\em arXiv preprint arXiv:1909.11929}, 2019.

\bibitem{krasikov2001nonnegative}
I.~Krasikov.
\newblock Nonnegative quadratic forms and bounds on orthogonal polynomials.
\newblock {\em Journal of Approximation Theory}, 111(1):31--49, 2001.

\bibitem{Csiszar}
I.~Csiszar and J.~K{\"o}rner.
\newblock {\em Information Theory: Coding Theorems for Discrete Memoryless
  Systems}.
\newblock Cambridge University Press, 2011.

\end{thebibliography}

\end{document}